\documentclass[12pt,a4paper]{article}

\usepackage{amssymb}
\usepackage{amsmath}
\usepackage{amsthm}
\usepackage{xypic}
\theoremstyle{plain}
\newtheorem{theorem}{Theorem}
\newtheorem{proposition}[theorem]{Proposition}
\newtheorem{corollary}[theorem]{Corollary}
\theoremstyle{remark}
\newtheorem{remark}[theorem]{Remark}
\theoremstyle{definition}
\newtheorem{definition}[theorem]{Definition}
\newtheorem{example}[theorem]{Example}
\newtheorem{lemma}[theorem]{Lemma}
\newtheorem{conjecture}[theorem]{Conjecture}

\def\varinjlim_#1{\lim\limits_{\longrightarrow\atop{#1}}}

\def\Hom{\mathop{\rm Hom}\nolimits}
\def\id{\mathop{\rm id}\nolimits}

\def\Z{\mathop{\rm Z}\nolimits}

\def\diag{\mathop{\rm diag}\nolimits}
\def\pt{\mathop{\rm pt}\nolimits}
\def\Ob{\mathop{\rm Ob}\nolimits}
\def\Mor{\mathop{\rm Mor}\nolimits}
\def\I{\mathop{\rm I}\nolimits}

\begin{document}

\author{A. V. Ershov}
\title{Formal groups over Hopf algebras}
\date{}
\maketitle
\section{The formal group, connected with FBSP}

The aim of this section is to define some generalization
of the notion of formal group. More precisely, we consider the
analog of formal groups with coefficients belonging to
a Hopf algebra. We also study some example of a formal group
over a Hopf algebra, which generalizes the formal group
of geometric cobordisms.

Recently some important connections
between the Landweber-Novikov algebra and the formal group
of geometric cobordisms were established (\cite{BBNY}).

Let $(H,\mu ,\eta ,\Delta ,\varepsilon, S)$ be a
(topological) Hopf algebra
over ring $R$ (where $\mu \colon H{\mathop{\widehat{\otimes}}\limits_R}H
\rightarrow H$ is the multiplication, $\eta \colon R\rightarrow H$
is the unit, $\Delta \colon H\rightarrow H
{\mathop{\widehat{\otimes}}\limits_R}H$
is the diagonal (comultiplication), $\varepsilon \colon H\rightarrow
R$ is the counit, and $S\colon H\rightarrow H$ is the antipode).

\begin{definition}
\label{fgoha}
A formal series ${\frak F}(x\otimes 1,1\otimes x)\in
H{\mathop{\widehat{\otimes}}\limits_R}H[[x\otimes 1,1\otimes x]]$
is called {\it a formal group over the Hopf algebra}
$(H,\mu ,\eta ,\Delta ,\varepsilon ,S)$
if the following conditions hold:
\begin{itemize}
\item[1)](associativity)
$$((\id_H\otimes \Delta){\frak F})
(x\otimes 1\otimes 1,1\otimes {\frak F}(x\otimes 1,1\otimes x))
=$$
$$((\Delta \otimes \id_H){\frak F})({\frak F}(x\otimes 1,
1\otimes x)\otimes 1,1\otimes 1\otimes x);$$
\item[2)](unit) $$((\id_H\otimes \varepsilon){\frak F})(x\otimes 1,
0)=x\otimes 1,$$ $$((\varepsilon \otimes \id_H){\frak F})(0,1\otimes x)
=1\otimes x;$$
\item[3)](inverse element) there exists the series
$\Theta(x)\in H[[x]]$ such that
$$((\mu \circ (\id_H\otimes S)){\frak F})
(x,\Theta(x))=0=((\mu \circ (S\otimes \id_H))
{\frak F})(\Theta(x),x).$$
\end{itemize}
If for a formal group
${\frak F}(x\otimes 1,1\otimes x)$ over
a commutative and cocommutative Hopf algebra $H$
the equality ${\frak F}(x\otimes 1,1\otimes x)=
{\frak F}(1\otimes x,x\otimes 1)$ holds, then it is called
commutative.
Below we shall deal only with the commutative case.
\end{definition}

\begin{remark}
Note that a formal group ${\frak F}(x\otimes 1,1\otimes x)$
over Hopf algebra $H$ over ring $R$
defines the formal group (in the usual sense) $F(x\otimes
1,1\otimes x)\in R[[x\otimes 1,1\otimes x]]$ over the ring $R$
in the following way. By $F(x\otimes 1,1\otimes x)$
denote the series
$((\varepsilon \otimes
\varepsilon){\frak F})(x\otimes 1,1\otimes x).$
If we identify $R{\mathop{\otimes}\limits_R}R$ and $R,$
we may assume that $F(x\otimes 1,1\otimes x)\in
R[[x\otimes 1,1\otimes x]].$
Note that for any coalgebra $H$ the diagram
\begin{equation}
\label{coalg}
\begin{array}{ccc}
H & \stackrel{\Delta}{\rightarrow} & H{\mathop{\widehat{\otimes}}
\limits_R}H\\
\scriptstyle{\varepsilon}
\downarrow & & \quad \downarrow
\scriptstyle{\varepsilon \otimes \varepsilon}\\
R & \stackrel{\cong}{\rightarrow} & R\otimes R\\
\end{array}
\end{equation}
is commutative. Using (\ref{coalg}) and
condition 1) of Definition \ref{fgoha}, we get $F(x\otimes 1\otimes 1,
1\otimes F(x\otimes 1,1\otimes x))=F(F(x\otimes 1,1\otimes x)\otimes 1,
1\otimes 1\otimes x).$
Similarly the conditions $F(x\otimes 1,0)=x\otimes 1$ and
$F(0,1\otimes x)=1\otimes x$
may be verified.
It is well known,
that the existence of the inverse element (in the case of usual
formal groups) follows from
the proved conditions.
However this may be deduced from the condition 3) of Definition
\ref{fgoha} in the standard way.
Moreover, the inverse element
$\theta(x)$ in the formal group $F(x\otimes 1,1\otimes x)$
is equal to $(\varepsilon(\Theta))(x).$

Therefore we may consider the formal group
${\frak F}(x\otimes 1,1\otimes x)$ over Hopf algebra $H$
as an extension of the usual formal group $F(x\otimes 1,1\otimes x)$
by the Hopf algebra $H.$
\end{remark}

\begin{remark}
By definition, put
$\widetilde{\Delta}(x)={\frak F}(x\otimes 1,1\otimes x),\;
\widetilde{\varepsilon}(x)=0,\; \widetilde{S}(x)=\Theta(x)$
and $\widetilde{\Delta}\mid_H=\Delta, \;
\widetilde{\varepsilon}\mid_H=\varepsilon, \;
\widetilde{S}\mid_H=S.$
We claim that $(H[[x]],\widetilde{\mu},\widetilde{\eta},
\widetilde{\Delta},\widetilde{\varepsilon},\widetilde{S})$
is the Hopf algebra
(here $\widetilde{\mu},\; \widetilde{\eta}$ are evidently
extensions of $\mu,\; \eta$).
Indeed, the commutativity of the diagram
\begin{equation}
\begin{array}{ccc}
H[[x]] & \stackrel{\widetilde{\Delta}}{\rightarrow} &
H[[x]]{\mathop{\widehat{\otimes}}\limits_R}H[[x]] \\
{\scriptstyle \widetilde{\Delta}}\downarrow \; \: &&
\qquad \qquad \downarrow {\scriptstyle \id_{H[[x]]}
\otimes \widetilde{\Delta}} \\
H[[x]]{\mathop{\widehat{\otimes}}\limits_R}H[[x]] &
\stackrel{\widetilde{\Delta}\otimes \id_{H[[x]]}}
{\rightarrow} & H[[x]]{\mathop{\widehat{\otimes}}\limits_R}
H[[x]]{\mathop{\widehat{\otimes}}\limits_R}H[[x]] \\
\end{array}
\end{equation}
follows from the equations
$$(\id_{H[[x]]}\otimes \widetilde{\Delta})
(\frak F(x\otimes 1,1\otimes x))=
((\id_H \otimes \Delta){\frak F})(x\otimes 1\otimes 1,
1\otimes {\frak F}(x\otimes 1,1\otimes x))=$$
$$((\Delta \otimes \id_H){\frak F})({\frak F}(x\otimes 1,1\otimes x)
\otimes 1,1\otimes 1\otimes x)=(\widetilde{\Delta}\otimes \id_{H[[x]]})
({\frak F}(x\otimes 1,1\otimes x)).$$
The commutativity of the diagram
$$
\diagram
R{\mathop{\otimes}\limits_R}H[[x]] &
\qquad H[[x]]{\mathop{\widehat{\otimes}}\limits_R}H[[x]]
\qquad \lto_{\widetilde{\varepsilon}\otimes \id_{H[[x]]}}
\rto^{\id_{H[[x]]} \otimes \widetilde{\varepsilon}}&
H[[x]]{\mathop{\otimes}\limits_R}R \\
& H[[x]]\ulto^\cong \uto_{\widetilde{\Delta}} \urto_\cong & \\
\enddiagram
$$
follows from the equations
$$((\id_{H[[x]]}\otimes \widetilde{\varepsilon})\circ \widetilde{\Delta})
(x)=((\id_H \otimes \varepsilon){\frak F})(x\otimes 1,1\otimes
\widetilde{\varepsilon}(x))=x\otimes 1,$$
$$((\widetilde{\varepsilon}\otimes \id_{H[[x]]})\circ \widetilde{\Delta})
(x)=((\varepsilon \otimes \id_H){\frak F})
(\widetilde{\varepsilon}(x)\otimes 1,
1\otimes x)=1\otimes x.$$
The axiom of antipode
$$(\widetilde{\mu}\circ (\id_{H[[x]]}\otimes \widetilde{S})
\circ \widetilde{\Delta})(x)=
(\widetilde{\mu}\circ (\widetilde{S}\otimes \id_{H[[x]]})
\circ \widetilde{\Delta})(x)=
(\widetilde{\eta}\circ \widetilde{\varepsilon})(x)=0$$
follows from the condition 3) of Definition \ref{fgoha}.
\end{remark}

\begin{remark}
We may rewrite the conditions 1), 2), 3)
of Definition \ref{fgoha} in terms
of series ${\frak F}(x\otimes 1,1\otimes x)$ in the next way.
Let $$\sum_{i,j\geq 0}
A_{i,j}x^i\otimes x^j=$$
$$\sum_{i,j\geq 0}
(\sum_ka_{i,j}^k
\otimes b_{i,j}^k)x^i\otimes x^j\in
H{\mathop{\widehat{\otimes}}\limits_R}H [[x\otimes 1,1\otimes
x]]$$ be the series ${\frak F}(x\otimes 1,1\otimes x).$ Then
the condition 1) is equivalent to the following equality:
$$\sum_{i,j\geq 0}
(\sum_ka_{i,j}^k\otimes
\Delta(b_{i,j}^k))x^i\otimes {\frak F}(x\otimes 1,1\otimes
x)^j=$$
$$\sum_{i,j\geq 0}
(\sum_k\Delta(a_{i,j}^k)
\otimes b_{i,j}^k){\frak F}(x\otimes 1,1\otimes x)^i\otimes
x^j$$
The condition 2) is equivalent to
$$\sum_ka_{i,0}^k
\varepsilon (b_{i,0}^k)=0,\quad {\rm if}\quad i\neq 1,\quad
\sum_ka_{1,0}^k\varepsilon (b_{1,0}^k)=1,$$
$$\sum_k\varepsilon (a_{0,j}^k)b_{0,j}^k=0,\quad
{\rm if}\quad j\neq 1,\quad \sum_k\varepsilon (a_{0,1}^k)
b_{0,1}^k=1.$$
The condition 3) also may be rewritten in terms of series.
\end{remark}

Let us consider some examples of defined objects.

\begin{example}\label{ttt}(Trivial extension)
Let $F(x\otimes 1,1\otimes x)$
be a formal group (in the usual sense) over a ring $R,$
and $(H,\mu ,\eta ,\Delta ,\varepsilon ,S)$
be a Hopf algebra over the same ring $R.$
Then ${\frak F}(x\otimes 1,1\otimes x)=
((\eta \otimes \eta )F)(x\otimes 1,1\otimes x)\in
H{\mathop{\widehat{\otimes}}\limits_R}H[[x\otimes 1,1\otimes x]]$
is the formal group over the Hopf algebra $H$
(recall that we identify $R{\mathop{\otimes}\limits_R}R$
and $R$).
\end{example}

\begin{example} Now we construct a nontrivial
extension ${\frak F}(x\otimes 1,1\otimes x)$ of the formal group of
geometric cobordisms
$F(x\otimes 1,1\otimes x)\in \Omega_U^*(\pt)[[x\otimes 1,1\otimes x]]$
by the Hopf algebra
$\Omega_U^*(Gr).$ For this let us consider the map
$\widehat{Gr}_{k,kl}\times \widehat{Gr}_{m,mn}
\stackrel{\widehat{\phi}_{kl,mn}}{\rightarrow}
\widehat{Gr}_{km,klmn},$ $(km,ln)=1.$
By $x\mid_{km,ln}$ denote the cobordism's class in
$\Omega_U^2(\widehat{Gr}_{km,klmn})$ such that
its restriction to every fiber of the bundle
\begin{equation}
\begin{array}{ccc}
\mathbb{C}P^{km-1} & \hookrightarrow & \widehat{Gr}_{km,klmn}\\
&& \downarrow \\
&& Gr_{km,klmn} \\
\end{array}
\end{equation}
is the standard generator in
$\Omega_U^2(\mathbb{C}P^{km-1}).$
Let $x\mid_{k,l}$ and $x\mid_{m,n}$ be analogously elements in
$\Omega_U^2(\widehat{Gr}_{k,l})$ and
$\Omega_U^2(\widehat{Gr}_{m,mn})$ respectively.
Then we obtain that $$\widehat{\phi}_{kl,mn}^*
(x\mid_{km,ln})=\sum_{0\leq i\leq k-1
\atop{0\leq j\leq m-1}}A_{i,j}
\mid_{kl,mn}(x\mid_{k,l})^i\otimes (x\mid_{m,n})^j,$$
where $A_{i,j}\mid_{kl,mn}\in \Omega_U^{2(1-i-j)}
(Gr_{k,kl}\times Gr_{m,mn}).$
Applying the functor of unitary cobordisms to the
following injective system
of the spaces and their maps
\begin{equation}
\begin{array}{ccc}
\widehat{Gr}_{p,pq}\times
\widehat{Gr}_{t,tu} &
\stackrel{\widehat{\phi}_{pq,tu}}{\rightarrow} &
\widehat{Gr}_{pt,pqtu}\\
\uparrow \quad && \uparrow \qquad \\
\widehat{Gr}_{k,kl}\times
\widehat{Gr}_{m,mn} &
\stackrel{\widehat{\phi}_{kl,mn}}{\rightarrow} &
\widehat{Gr}_{km,klmn}.\\
\end{array}
\end{equation}
(under the conditions $k\mid p,\quad l\mid q,\quad
m\mid t,\quad n\mid u,\quad$and
$(pt,qu)=1$),
we obtain the formal series $${\frak F}(x\otimes 1,1\otimes x)=$$
$$\sum_{i,j\geq 0}A_{i,j}x^i\otimes x^j\quad \in \Omega_U^*
(Gr){\mathop{\widehat{\otimes}}
\limits_{\Omega_U^*(\pt)}}\Omega_U^*(Gr)[[x\otimes 1,
1\otimes x]]$$
such that $i^*_{kl}A_{i,j}=A_{i,j}\mid_{k,l}$ for injection
$i_{kl}\colon Gr_{k,kl}\hookrightarrow Gr$ (for every pair $\{
k,l\}$ such that $(k,l)=1$).

By $R$ and $H$ denote the ring
$\Omega_U^*(\pt)$ and the Hopf algebra $\Omega_U^*(Gr)$ (over
the ring $\Omega_U^*(\pt)$) respectively
(recall that we consider the space $Gr$ with
the $H$-group structure, induced by the multiplication of
FBSP).

\begin{proposition}
\label{gfggc}
The series
${\frak F}(x\otimes 1,1\otimes x)$
is the formal group over the Hopf algebra $H.$
\end{proposition}
{\raggedright {\it Proof}.}
To prove $((\id_H\otimes \Delta){\frak F})(x\otimes 1\otimes 1,1\otimes
{\frak F}(x\otimes 1,1\otimes x))=((\Delta \otimes \id_H){\frak F})
({\frak F}(x\otimes 1,1\otimes x)\otimes 1,1\otimes 1\otimes x),$
we need the following commutative diagram ($(kmt,lnu)=1$):
\begin{equation}
\begin{array}{ccc}
\widehat{Gr}_{k,kl}\times \widehat{Gr}_{m,mn}
\times \widehat{Gr}_{t,tu} & \rightarrow &
\widehat{Gr}_{k,kl}\times
\widehat{Gr}_{mt,mntu}\\
\downarrow && \downarrow \\
\widehat{Gr}_{km,klmn}\times \widehat{Gr}_{t,tu} &
\rightarrow & \widehat{Gr}_{kmt,klmntu}. \\
\end{array}
\end{equation}

To prove $((\id_H\otimes \varepsilon){\frak F})(x\otimes 1,0)=
x\otimes 1,$
we need the following commutative diagram ($(km,ln)=1$):
\begin{equation}
\begin{array}{ccc}
\widehat{Gr}_{k,kl}\times \widehat{Gr}_{m,mn} & \rightarrow &
\widehat{Gr}_{km,klmn} \\
\uparrow && \uparrow \\
\widehat{Gr}_{k,kl}\times \mathbb{C}P^{m-1} &
\leftarrow & \widehat{Gr}_{k,kl}\times \{ \pt \}, \\
\end{array}
\end{equation}
where right-hand vertical arrow is the standard inclusion.

To prove $((\mu \circ (\id_H\otimes S)){\frak F})(x,\Theta
(x))=0,$ let us construct the fiber map $\widehat{\nu}\colon
\widehat{Gr}\rightarrow \widehat{Gr}$ such that the
following two conditions are satisfied:
\begin{itemize}
\item[1)] the restriction of $\widehat{\nu}$
to any fiber ($\cong \mathbb{C}P^\infty$)
is the inversion in the $H$-group $\mathbb{C}P^\infty$;
\item[2)] $\widehat{\nu}$ covers the $\nu \colon Gr\rightarrow Gr$
(where $\nu$ is the inversion in the $H$-group $Gr$).
\end{itemize}

Let us remember that
${\cal P}^{k-1}{\mathop{\times}\limits_{Gr_{k,kl}}}{\cal Q}^{l-1}$
is the canonical FBSP over $Gr_{k,kl}$ and we have denoted
by $\widehat{Gr}_{k,kl}$ the bundle space ${\cal P}^{k-1}.$
Let $\widehat{Gr}'_{k,kl}$ ($\widehat{Gr}'$) be the bundle space of the
''second half'' ${\cal Q}^{l-1}$ of the canonical FBSP over
$Gr_{k,kl}$ ($\varinjlim_{(k,l)=1}\widehat{Gr}'_{k,kl}$
respectively).

First note that there exists the fiber isomorphism
$\widehat{\nu}'_{k,l}\colon \widehat{Gr}_{k,kl}
\rightarrow \widehat{Gr}'_{l,lk}$
that covers the inverse map $\nu_{k,l}\colon Gr_{k,kl}\rightarrow
Gr_{l,lk}$ (in other words, the map $\nu_{k,l}$ takes each
subalgebra $A_k\cong M_k(\mathbb{C})$ in the $M_{kl}(\mathbb{C})$
to its centralizer $\Z_{M_{kl}(\mathbb{C})}(A_k)\cong M_l(\mathbb{C})$
in the $M_{kl}(\mathbb{C})$). Let $c_{l,k}\colon
\widehat{Gr}'_{l,lk}\rightarrow \widehat{Gr}'_{l,lk}$
be the fiber map such that the following two conditions
are satisfied:
\begin{itemize}
\item[1)] $c_{l,k}$ covers the identity mapping of the base $Gr_{l,lk};$
\item[2)] the restriction of $c_{l,k}$ to any fiber
$\cong \mathbb{C}P^{k-1}$ is the complex conjugation.
\end{itemize}
Let $\widehat{\nu}_{k,l}\colon \widehat{Gr}_{k,kl}\rightarrow
\widehat{Gr}'_{l,lk}$ be the composition
$c_{l,k}\circ \widehat{\nu}'_{k,l}.$
It is easy to prove that the map
$\widehat{\nu}=\varinjlim_{(k,l)=1}\widehat{\nu}_{k,l}
\colon \varinjlim_{(k,l)=1}\widehat{Gr}_{k,kl}
\rightarrow \varinjlim_{(l,k)=1}\widehat{Gr}'_{l,lk}$
is required. In particular, there exists the fiber isomorphism
between $\widehat{Gr}$ and $\widehat{Gr}'.$

The map $\widehat{\nu}$ defines (by the same way, as
$\widehat{\phi}$ in the beginning of the example) the formal series
$\Theta(x)\in H[[x]]$ (note that
$\varepsilon(\Theta)(x)=
\theta(x),$ where $\theta(x)\in R[[x]]$ is the inverse element in
the group of geometric cobordisms).

Now we claim that
$((\mu \circ (\id_H\otimes S)){\frak F})(x,\Theta(x))=0.$
Indeed, this follows from the next commutative diagram:
\begin{equation}
\begin{array}{ccccccc}
\widehat{Gr} & \stackrel{\diag}{\rightarrow} &
\widehat{Gr}\times \widehat{Gr} &
\stackrel{\id \times \widehat{\nu}}{\rightarrow} &
\widehat{Gr}\times \widehat{Gr} & \stackrel{\widehat{\phi}}
{\rightarrow} & \widehat{Gr} \\
\downarrow && \downarrow && \downarrow && \downarrow \\
Gr & \stackrel{\diag}{\rightarrow} & Gr\times Gr &
\stackrel{\id \times \nu}{\rightarrow} & Gr\times Gr &
\stackrel{\phi}{\rightarrow} & Gr \\
\end{array}
\end{equation}
(we see that the composition
$\widehat{\phi} \circ (\id \times \widehat{\nu})\circ \diag$
is homotopic (in class of fiber homotopies) to the
map $\widehat{Gr}\rightarrow \pt \in \widehat{Gr}$).
$\square$

\smallskip

Let $\kappa \colon
\widetilde{Gr}\rightarrow \mathbb{C}P^\infty$
be the direct limit of the fiber maps
\begin{equation}
\begin{array}{ccc}
\widetilde{Gr}_{k,kl} & \stackrel{\kappa_{k,l}}{\rightarrow} &
\mathbb{C}P^{kl-1} \\
\downarrow && \downarrow \\
Gr_{k,kl} & \rightarrow & \pt. \\
\end{array}
\end{equation}
It defines (in the same way, as $\widehat{\phi}$ and
$\widehat{\nu}$ above) the formal
series
$${\frak G}(x,y)=
\sum_{i,j\geq 0}B_{i,j}x^iy^j\; \in H[[x,y]].$$

\begin{proposition}
${\frak G}(x,y)=
((\mu \circ (\id_H\otimes S)){\frak F})(x,y),$
i. e. $$B_{i,j}=\sum_ka_{i,j}^kS(b_{i,j}^k).$$
\end{proposition}
{\raggedright {\it Proof}.}
Recall that in the proof of Proposition \ref{gfggc}
the fiber maps $\widehat{\nu}'_{k,l}\colon \widehat{Gr}_{k,kl}
\rightarrow \widehat{Gr}'_{l,lk}$ were defined.
By $\widehat{\nu}'$ denote the direct limit $\varinjlim_{(k,l)=1}
\widehat{\nu}'_{k,l}\colon \widehat{Gr}\rightarrow \widehat{Gr}.$
Note that $\widehat{\nu}'$ covers the inversion
$\nu \colon Gr\rightarrow Gr$ in the $H$-group $Gr.$

Now the proof follows from the next composition of the bundle maps:
\begin{equation}
\begin{array}{ccccccc}
\widetilde{Gr} & \rightarrow &
\widehat{Gr}\times \widehat{Gr} & \stackrel{\id \times \widehat{\nu}'}
{\rightarrow} &
\widehat{Gr}\times
\widehat{Gr} & \stackrel{\widehat{\phi}}{\rightarrow} & \widehat{Gr} \\
\downarrow && \downarrow && \downarrow && \downarrow \\
Gr & \stackrel{diag}{\rightarrow} & Gr\times Gr & \stackrel{\id \times
\nu}{\rightarrow} & Gr\times Gr & \stackrel{\phi}{\rightarrow} & Gr. \\
\end{array}
\end{equation}
We see that the upper composition in fact is the map $\widetilde{Gr}
\rightarrow \mathbb{C}P^\infty$ and it coincides with the map $\kappa.$
Let $y$ be $\widehat{\nu}'{^*}(x).$ The upper composition gives
$x\mapsto {\frak F}(x\otimes 1,1\otimes x)\mapsto
((\id_H\otimes S){\frak F})(x\otimes 1,1\otimes y)\mapsto
((\mu \circ (\id_H\otimes S)){\frak F})(x\otimes 1,1\otimes y).$
Without loss of sense we may write $x$ and $y$ instead of
$x\otimes 1$ and $1\otimes y$ respectively. $\square$

The series ${\frak G}(x,y)$ has the following
interesting property.
\begin{proposition}
$$(\Delta{\frak G})({\frak F}(x\otimes 1,1\otimes x),
((S\otimes S){\frak F})(y\otimes 1,1\otimes y))=$$
$$F({\frak G}(x,y)\otimes 1,1\otimes{\frak G}(x,y)),$$
where $F(x,y)\; \in R[[x,y]]$ is the formal group
of geometric cobordisms.
\end{proposition}
{\raggedright {\it Proof}.}
We give two variants of the proof.

1).''Topological proof'' follows from the commutative diagram
\begin{equation}
\begin{array}{ccc}
\mathbb{C}P^{kl-1}\times \mathbb{C}P^{mn-1} & \rightarrow &
\mathbb{C}P^{klmn-1} \\
\uparrow \; \: && \uparrow \\
\widetilde{Gr}_{k,kl}\times \widetilde{Gr}_{m,mn} & \rightarrow &
\widetilde{Gr}_{km,klmn} \\
\end{array}
\end{equation}
($(km,ln)=1$) combining with the decomposition of the map $\kappa,$
which was obtained in previous proof.

2). By $\widetilde{S}'$ denote the homomorphism
$\widehat{\nu}'{^*}
\colon H[[x]]\rightarrow H[[y]]$ (recall that $\widehat{\nu}'{^*}\mid_H=S
\colon H\rightarrow H,$ where $S$ is the antipode).
Let us consider the following composition of homomorphisms
of Hopf algebras:
$$H[[x]]\stackrel{\widetilde{\Delta}}{\rightarrow}
H[[x]]{\mathop{\widehat{\otimes}}\limits_R}H[[x]]
\stackrel{\id \otimes \widetilde{S}'}{\rightarrow}
H[[x]]{\mathop{\widehat{\otimes}}\limits_R}H[[y]]
\stackrel{(\mu)}{\rightarrow}H[[x,y]],$$
where $(\mu)$ is the homomorphism, induced by multiplication
$$\mu \colon H{\mathop{\widehat{\otimes}}\limits_R}H\rightarrow H.$$
It follows from the axiom of antipode $$\mu \circ (\id_H\otimes S)\circ
\Delta =\eta \circ \varepsilon$$ that $$(\mu)\circ (\id_{H[[x]]}\otimes
\widetilde{S}')\circ \widetilde{\Delta}\mid_H=\eta \circ \varepsilon.$$
Hence there exists the homomorphism of Hopf algebras
$$(\eta)\colon R[[x]]\rightarrow H[[x,y]]$$
such that the following diagram
$$
\diagram
H[[x]] \rto^{{\widetilde{\Delta}}\quad}
\drto_{(\varepsilon)} & H[[x]]{\mathop{\widehat{\otimes}}\limits_R}H[[x]]
\rto^{\id \otimes \widetilde{S}'} &
H[[x]]{\mathop{\widehat{\otimes}}\limits_R}H[[y]] \dto^{(\mu)} \\
& R[[x]] \rto^{(\eta)\qquad} & H[[x,y]] \\
\enddiagram
$$
is commutative (here $(\varepsilon)$ is the homomorphism,
induced by $\varepsilon$). Note that $$(\eta)(x)=
{\frak G}(x,y).$$
This completes the proof that $$\Delta_{R[[x]]}(x)=F(x\otimes 1,
1\otimes x),$$ where $F(x\otimes 1,1\otimes x)\in R[[x\otimes
1,1\otimes x]]$ is the formal group of geometric cobordisms. $\square$

It is very important that we consider the maps
$\widehat{\phi},\; \widehat{\nu},$ and $\kappa$ as {\it fiber} maps
in this example. Otherwise instead of
${\frak F}(x\otimes 1,1\otimes x)$ we obtain the usual formal group
of geometric cobordisms because the $H$-space $\widehat{Gr}$
is isomorphic to the $H$-space $BSU_\otimes \times \mathbb{C}P^\infty.$
\end{example}

It is well known (\cite{Quillen}), that the formal group
of geometric cobordisms is the universal
formal group.
\begin{conjecture}
\label{conj}
The formal group ${\frak F}(x\otimes 1,1\otimes x)$
is the universal object in the category of formal groups
over a (topological) Hopf algebras.
\end{conjecture}

Let $R'$ be a ring and $F'(x\otimes 1,1\otimes x)$
be a formal group over $R'.$ Note that we may
consider the $R'$ as the Hopf algebra over $R'$
with respect to the $\Delta_{R'}\colon R'\cong R'{\mathop{\otimes}
\limits_{R'}}R',$ $\eta_{R'}=\varepsilon_{R'}=S_{R'}=\id_{R'}
\colon R'\rightarrow R'.$
If $\chi \colon H\rightarrow R'$ is a homomorphism of the
Hopf algebras from $(H,\mu,\eta,\Delta,\varepsilon,S)$ to
$(R',\mu_{R'},\eta_{R'},\Delta_{R'},\varepsilon_{R'},S_{R'}),$
then $\chi = (\chi \circ \eta)\circ \varepsilon =\chi \mid _R
\circ \varepsilon.$ Hence there exists the natural bijection
$\Hom_{Hopf\: alg.}(H,R')\leftrightarrow \Hom_{Ring}(R,R').$
Therefore the Conjecture
implies the universal property of the
formal group of geometric cobordisms.

\smallskip

\section{Extensions
of the formal group of geometric cobordisms,
generated by ${\frak F}(x\otimes 1,1\otimes x)$}

Now we construct the denumerable set of extensions
of the formal group of geometric cobordisms
$F(x\otimes 1,1\otimes x)$
by the Hopf algebra $H=\Omega_U^*(Gr).$

Let $F_i(x\otimes 1,1\otimes x),\; i=1,2$
be formal groups over ring $R.$
Recall the following definition.
\begin{definition}
A homomorphism of formal groups
$\varphi \colon F_1\rightarrow F_2$
is a formal series $\varphi(x)\in R[[x]]$ such that
$\varphi(F_1(x\otimes 1,1\otimes x))=F_2(\varphi(x)\otimes 1,
1\otimes \varphi(x)).$
\end{definition}
Let $H$ be a Hopf algebra over ring $R$ with
diagonal $\Delta;$
let ${\frak F}_i(x\otimes 1,1\otimes x),\; i=1,2$ be
formal groups over $H.$
\begin{definition}
A homomorphism of formal groups over Hopf algebra $H$
$\Phi \colon {\frak F}_1\rightarrow {\frak F}_2$
is a formal series $\Phi(x)\in H[[x]]$ such that
$(\Delta \Phi)({\frak F}_1(x\otimes 1,1\otimes x))={\frak F}_2(\Phi(x)
\otimes 1,1\otimes\Phi(x)).$
\end{definition}
Note that $\varepsilon(\Phi)\colon
(\varepsilon \otimes \varepsilon)({\frak F}_1)
\rightarrow (\varepsilon \otimes \varepsilon)({\frak F}_2)$
is the homomorphism of the formal groups over the ring $R$
(where $\varepsilon$ is the counit of the Hopf algebra $H$).
We say that the homomorphism $\Phi$ covers the homomorphism
$\varepsilon(\Phi).$

Let $R$ be the ring $\Omega_U^*(\pt);$ let
$F(x\otimes 1,1\otimes x)\in R[[
x\otimes 1,1\otimes x]]$
be the formal group of geometric cobordisms.
Let $H$ be the Hopf algebra $\Omega_U^*(Gr).$
By definition, put
$\varphi^{(1)}(x)=x,\; \varphi^{(-1)}(x)=\theta(x)$
and $\varphi^{(n)}(x)=F(x,\varphi^{(n-1)}(x)),$
where $\theta(x)\in R[[x]]$ is the inverse element in $F.$
Clearly, that $\varphi^{(n)}\colon F\rightarrow F$
is the homomorphism for every $n\in \mathbb{Z}.$
Power systems were considered by S.~P.~Novikov
and V.~M.~Buchstaber in \cite{novbuch}.

Below for any $n\in \mathbb{Z}$ we construct the
extension ${\frak F}^{(n)}(x\otimes 1,1\otimes x)$
of $F(x\otimes 1,1\otimes x)$ by $H$
and the homomorphism $\Phi^{(n)}\colon {\frak F}\rightarrow
{\frak F}^{(n)}$ such that
\begin{itemize}
\item[(i)] ${\frak F}^{(1)}={\frak F};$
\item[(ii)] $\varepsilon(\Phi^{(n)})=\varphi^{(n)}.$
\end{itemize}

Let $n$ be a positive integer.
Let us take the product of the FBSP $\widetilde{Gr}_{k,kl}$
(over $Gr_{k,kl}$) with itself $n$ times. It is the FBSP
over $Gr_{k,kl}$ with a fiber $\mathbb{C}P^{k^n-1}\times
\mathbb{C}P^{l^n-1}.$ By $\widetilde{Gr}_{k,kl}^{(n)}$
denote the obtained FBSP.
Let $\widehat{Gr}_{k,kl}^{(n)}$ be the corresponding
bundle over $Gr_{k,kl}$ with fiber $\mathbb{C}P^{k^n-1}.$
Let $\widehat{Gr}^{(n)}=\varinjlim_{(k,l)=1}\widehat{Gr}_{k,kl}^{(n)}.$
We have the evident fiber maps $\widehat{Gr}_{k,kl}\rightarrow
\widehat{Gr}_{k,kl}^{(n)},\quad \lambda^{(n)} \colon \widehat{Gr}\rightarrow
\widehat{Gr}^{(n)}$ and
the following commutative diagrams ($(km,ln)=1$):
\begin{equation}
\begin{array}{ccc}
\widehat{Gr}_{km,klmn} & \rightarrow & \widehat{Gr}_{km,klmn}^{(n)} \\
\uparrow && \uparrow \\
\widehat{Gr}_{k,kl}\times \widehat{Gr}_{m,mn} &
\rightarrow & \widehat{Gr}_{k,kl}^{(n)}\times
\widehat{Gr}_{m,mn}^{(n)},\\
\end{array}
\end{equation}
\begin{equation}
\label{2}
\begin{array}{ccc}
\widehat{Gr} & \stackrel{\lambda^{(n)}}{\rightarrow} & \widehat{Gr}^{(n)}\\
\scriptstyle{\widehat{\phi}}\uparrow && \uparrow
\scriptstyle{\widehat{\phi}^{(n)}} \\
\widehat{Gr}\times \widehat{Gr} & \stackrel{\lambda^{(n)}
\times \lambda^{(n)}}{\rightarrow} &
\widehat{Gr}^{(n)}\times \widehat{Gr}^{(n)}.\\
\end{array}
\end{equation}
By $x$ denote the class of cobordisms in
$\Omega_U^{2}(\widehat{Gr}^{(n)})$
such that its restriction to any fiber $\cong \mathbb{C}P^\infty$
is the generator
$x\mid _{\mathbb{C}P^\infty}\in \Omega_U^2(\mathbb{C}P^\infty)$.
Let $\Phi^{(n)}(x)\in H[[x]]$ be the series, defined by
the fiber map $\lambda^{(n)}.$
Let
$$
{\frak F}^{(n)}(x\otimes 1,1\otimes x)\in H{\mathop{\widehat{\otimes}}
\limits_R}H[[x\otimes 1,1\otimes x]]
$$
be the series,
corresponds to the fiber map $\widehat{Gr}^{(n)}\times
\widehat{Gr}^{(n)}\stackrel{\widehat{\phi}^{(n)}}
{\rightarrow}\widehat{Gr}^{(n)};$
note that $Gr^{(n)}$ is the $H$-group
with the multiplication $\widehat{\phi}^{(n)}$).
Clearly, that ${\frak F}^{(n)}(x\otimes 1,1\otimes x)$
is an extension of $F(x\otimes 1,1\otimes x)$ by $H$
(in particular, it is the formal group over Hopf algebra $H$).
Note that $\lambda^{(n)}$ covers the identity map
of the base $Gr.$ It follows from diagram (\ref{2}) that
$$
(\Delta \Phi^{(n)})({\frak F}(x\otimes 1,1\otimes x))=
{\frak F}^{(n)}(\Phi^{(n)}(x)\otimes 1,1\otimes \Phi^{(n)}(x)).
$$
It is clear that $\varepsilon(\Phi^{(n)})(x)=\varphi^{(n)}(x).$

For $n=0$ let $\widehat{Gr}^{(0)}=Gr\times \mathbb{C}P^{\infty}$
and let $\lambda^{(0)}$ be the composition
$$\widehat{Gr}\rightarrow \pt \rightarrow \widehat{Gr}^{(0)}.$$
It defines the series ${\frak F}^{(0)}=F$ and $\Phi^{(0)}=0.$

Let $\lambda^{(-1)}$ be the fiber map $\widehat{Gr}\rightarrow
\widehat{Gr}^{(-1)}=\widehat{Gr}$ such that the following conditions
hold:
\begin{itemize}
\item[(i)] the restriction of $\lambda^{(-1)}$ to any fiber
is the inversion in the $H$-group $\mathbb{C}P^{\infty}$
(i. e. the complex conjugation);
\item[(ii)] $\lambda^{(-1)}$ covers the map $\nu \colon Gr \rightarrow Gr,$
where $\nu$ is the inversion in the $H$-group $Gr.$
\end{itemize}
Let $\Phi^{(-1)}(x)\in H[[x]]$ be the series,
defined by $\lambda^{(-1)}.$
Trivially, that $\varepsilon(\Phi^{(-1)})(x)=\theta(x).$
Note that the $\lambda^{(-1)}$ coincides with $\widehat{\nu}.$
Consequently, $\Phi^{(-1)}=\Theta(x).$
Now we can define ${\frak F}^{(n)}$ and $\Phi^{(n)}$
for negative integer $n$ by the obvious way.

By $S$ denote the antipode
of the Hopf algebra $H.$ Let $\mu$ be the multiplication
in the Hopf algebra $H.$
By definition, put $(1)=\id_H,\; (-1)=S\colon H\rightarrow H$ and
$(n)=\mu \circ ((n-1)\otimes (1))\circ \Delta \colon H\rightarrow H$
(in particular, $(0)=\eta \circ \varepsilon \colon H\rightarrow H,$
where $\eta$ is the unit in $H$).
\begin{proposition}
${\frak F}^{(n)}(x\otimes 1,1\otimes x)=
(((n)\otimes(n)){\frak F})(x\otimes 1,1\otimes x)$
for any $n\in \mathbb{Z}.$
\end{proposition}
{\raggedright {\it Proof}.}
By $\phi \colon Gr\times Gr\rightarrow Gr$ denote the multiplication
in the $H$-space $Gr.$
Suppose $n$ a positive integer.
By definition, put $\phi{(1)}=\id_{Gr},\;
\phi{(n)}=\phi \circ (\phi{(n-1)} \times \id_{Gr}),$
and $\diag{(n)}=(\diag{(n-1)}\times \id_{Gr})\circ \diag ,$
where $\diag(1)=\id_{Gr},\;
\diag =\diag{(2)}\colon Gr\rightarrow Gr\times Gr.$
Note that the composition $\phi{(n)}\circ \diag{(n)}\colon Gr
\rightarrow Gr$ induces the homomorphism $(n)\colon H\rightarrow H.$

Let us consider the classifying
map $\alpha{(n)}\colon Gr\rightarrow Gr$ for the bundle $\widehat{Gr}^{(n)}$
over $Gr.$ We have the following commutative diagram:
$$
\diagram
\mathbb{C}P^\infty \drto \rto^= & \mathbb{C}P^\infty \drto \\
& \quad \widehat{Gr}^{(n)} \dto \rto^{\widehat{\alpha}(n)} &
\quad \widehat{Gr} \dto \\
& Gr \rto^{\alpha{(n)}} & Gr \\
\enddiagram
$$
It is easy to prove that $\alpha{(n)}=\phi{(n)}\circ \diag{(n)}.$
Hence $\alpha{(n)}^*=(n)\colon H\rightarrow H.$
Note that the following diagram
\begin{equation}
\begin{array}{ccc}
\widehat{Gr}^{(n)}\times \widehat{Gr}^{(n)} &
\stackrel{\widehat{\alpha}{(n)}
\times \widehat{\alpha}{(n)}}{\rightarrow} &
\widehat{Gr}\times \widehat{Gr}\\
\scriptstyle{\widehat{\phi}^{(n)}}\downarrow \qquad &&
\quad \downarrow \scriptstyle{\widehat{\phi}} \\
\widehat{Gr}^{(n)} & \stackrel{\widehat{\alpha}{(n)}}
{\rightarrow} & \widehat{Gr}\\
\end{array}
\end{equation}
is commutative. This completes the proof for positive $n.$
For negative $n$ proof is similar. $\square$

\smallskip

We can define the structure of group on the set $\{ {\frak
F}^{(n)};\; n\in \mathbb{Z}\}$ in the following way.
Recall that for any Hopf algebra $H$ the triple $(\Hom_{Alg. Hopf}(H,H),
\star, \eta \circ \varepsilon)$
is the algebra with respect to the convolution $f\star g=
\mu \circ (f\otimes g)\circ \Delta \colon H\rightarrow H.$
It follows from the previous Proposition that the formal group
${\frak F}^{(n)}$ corresponds to the homomorphism $(n)\colon H
\rightarrow H$ (see Conjecture \ref{conj}).
Clearly, that $(m)\star(n)=(m+n)$ for any $m,n\in \mathbb{Z}.$

\section{Logarithms of formal groups \\ over Hopf algebras}

In this section by $(H,\mu,\eta,\Delta,\varepsilon,S)$
denote a commutative Hopf algebra over ring $R$ without torsion
and by ${\frak F}(x\otimes 1,1\otimes x)$ denote a
formal group over Hopf algebra $H.$
By $H_\mathbb{Q}$ denote the Hopf algebra
$H{\mathop{\otimes}\limits_\mathbb{Z}}\mathbb{Q}$
over ring
$R_\mathbb{Q}=R{\mathop{\otimes}\limits_\mathbb{Z}}\mathbb{Q}.$
We shall write $\mu,\eta,\ldots$ instead of $\mu_\mathbb{Q},
\eta_\mathbb{Q},\ldots.$

The aim of this section is to prove the following result.
\begin{proposition}
For any commutative formal group ${\frak F}(x\otimes 1,1\otimes x),$
which is considered as a formal group over $H_\mathbb{Q},$
there exists a homomorphism  to
a formal group of the form ${\frak c}+x\otimes 1+1\otimes x,$
where $\frak c\in H_\mathbb{Q}{\mathop{\widehat{\otimes}}
\limits_{R_\mathbb{Q}}}H_\mathbb{Q}$
such that $(\id \otimes \varepsilon){\frak c}=0=
(\varepsilon \otimes \id){\frak c}.$
\end{proposition}
We recall that the notion of a homomorphism of formal
groups over Hopf algebra was given in previous Section.
Below we shall use notations of Section~1.

To prove the Proposition, we need the following Lemma.
\begin{lemma}
A symmetric series of the form ${\frak c}+x\otimes 1+1\otimes x
\in H_\mathbb{Q}{\mathop{\widehat{\otimes}}
\limits_{R_\mathbb{Q}}}H_\mathbb{Q}[[x\otimes 1,1\otimes x]]$
is a formal group over the Hopf algebra $H_\mathbb{Q}$
if and only if the following two conditions hold:
\begin{itemize}
\item[(i)] $(\id \otimes \Delta){\frak c}+1\otimes
{\frak c}-(\Delta \otimes \id){\frak c}-{\frak c}\otimes 1=0;$
\item[(ii)] $(\id \otimes \varepsilon){\frak c}=0=
(\varepsilon \otimes \id){\frak c}.$
\end{itemize}
\end{lemma}
{\raggedright(Note} that the condition (i) means, that ${\frak c}$ is a
2-cocykle in the cobar complex of the Hopf algebra
$H_\mathbb{Q}.$)\\
{\raggedright {\it Proof of the Lemma}.}\;
The conditions (i) and (ii) are equivalent to the associativity
axiom and to the unit axiom for formal groups respectively.
Let us show that the series $\Theta(x)=-(\mu\circ (\id \otimes S))
{\frak c}-x$ is the inverse element. Indeed,
$$({\mu}\circ (\id \otimes S)){\frak c}+
x+\Theta(x)=
(\mu \circ (\id \otimes S)){\frak c}+x-
(\mu \circ (\id \otimes S)){\frak c}-x=0.$$
The symmetric condition follows from the equality
$({\mu}\circ (\id \otimes S)){\frak
c}=({\mu}\circ (S\otimes \id)){\frak c}.\quad \square$

{\raggedright {\it Proof of the Proposition}.}\;
By definition, put
$$\widetilde{\omega}(x)=(\id \otimes \widetilde{\varepsilon})
\frac{\partial{\frak F}(x,z)}{\partial z}\in H[[x]]$$
(here $\widetilde{\varepsilon}\colon H[[z]]\rightarrow R$
is the map such that $\widetilde{\varepsilon}\mid_H=\varepsilon
\colon H\rightarrow R,\; \widetilde{\varepsilon}(z)=0$).
Recall that $\widetilde{\Delta}\colon
H[[x]]\rightarrow H[[x]]{\mathop{\widehat{\otimes}}\limits_R}H[[x]]=
H{\mathop{\widehat{\otimes}}\limits_R}H[[x\otimes 1,1\otimes x]]$
is the map such that $\widetilde{\Delta}\mid_H=\Delta,\;
\widetilde{\Delta}(x)={\frak F}(x\otimes 1,1\otimes x).$
We have
$$
(\Delta \widetilde{\omega})({\frak F}(x\otimes 1,1\otimes x))=
\widetilde{\Delta}(\widetilde{\omega}(x))=(\widetilde{\Delta}
\circ (\id \otimes \widetilde{\varepsilon})\circ \frac{\partial}
{\partial z})({\frak F}(x,z))=
$$
$$
((\id \otimes \id \otimes \widetilde{\varepsilon})\circ
(\widetilde{\Delta}\otimes \id)\circ \frac{\partial}{\partial z})
({\frak F}(x,z))=
$$
$$
((\id \otimes \id \otimes \widetilde{\varepsilon})\circ
\frac{\partial}{\partial z}\circ (\widetilde{\Delta}\otimes
\id))({\frak F}(x,z))=
$$
$$
\left((\id \otimes \id \otimes \widetilde{\varepsilon})\circ
\frac{\partial}{\partial z}\right)((\Delta \otimes \id){\frak F})({\frak
F}(x\otimes 1,1\otimes x),z)=
$$
$$
\left((\id \otimes \id \otimes \widetilde{\varepsilon})\circ
\frac{\partial}{\partial z}\right)((\id \otimes \Delta){\frak F})(x
\otimes 1,{\frak F}(1\otimes x,z))=
$$
$$
(\id \otimes \id \otimes \widetilde{\varepsilon})\left(\frac{\partial
((\id \otimes \Delta){\frak F})(x\otimes 1,{\frak F}(1\otimes x,
z))}{\partial {\frak F}(1\otimes x,z)}\right)\cdot 1\otimes\left((\id
\otimes \widetilde{\varepsilon})\frac{\partial{\frak F}
(1\otimes x,z)}{\partial z}\right)=
$$
$$
\frac{\partial{\frak F}(x\otimes 1,1\otimes x)}
{\partial (1\otimes x)}\cdot(1\otimes \widetilde{\omega})(1\otimes x).
$$
Therefore, we have
\begin{equation}
\label{1}
(\Delta \widetilde{\omega})({\frak F}(x\otimes 1,1\otimes x))=
\frac{\partial{\frak F}(x\otimes 1,1\otimes x)}
{\partial (1\otimes x)}\cdot(1\otimes \widetilde{\omega})(1\otimes x).
\end{equation}
If $${\frak F}(x,z)=\sum_{i,j\geq 0}A_{i,j}x^iz^j\quad (A_{i,j}\in
H{\mathop{\widehat{\otimes}}\limits_R}H),$$
then
$$\widetilde{\omega}(x)=(\id \otimes \widetilde{\varepsilon})
\sum_{i,j}A_{i,j}x^ijz^{j-1}=(\id \otimes \varepsilon)A_{0,1}+
\sum_{i\geq 1}((\id \otimes \varepsilon)A_{i,1})x^i,
$$
where $(\varepsilon \circ(\id \otimes \varepsilon))A_{0,1}=1\neq 0.$
Therefore
$$\frac{1}{\widetilde{\omega}(x)}\in H[[x]]\quad \hbox{and}
$$
$$
\widetilde{\Delta}\left(\frac{1}{\widetilde{\omega}(x)}\right)=
\frac{1}{\widetilde{\Delta}(\widetilde{\omega}(x))}=
\frac{1}{(\Delta\widetilde{\omega})({\frak F}(x\otimes 1,1\otimes x))}
\in H{\mathop{\widehat{\otimes}}\limits_R}H[[x\otimes 1,1\otimes x]].
$$
Therefore (\ref{1}) may be rewritten in the form
\begin{equation}
\label{2}
\frac{d(1\otimes x)}{(1\otimes
\widetilde{\omega})(1\otimes x)}=\frac{d{\frak F}
(x\otimes 1,1\otimes x)}{(\Delta\widetilde{\omega})({\frak F}
(x\otimes 1,1\otimes x))}.
\end{equation}
It is clear that
$$
\frac{1}{\widetilde{\omega}(x)}=b_0+
b_1x+\ldots ,
$$
where $b_i\in H,\; \varepsilon(b_0)=1.$
By ${\frak g}(x)$ denote the series
$$\int_o^x\frac{dt}{\widetilde{\omega}(t)}
\in H_\mathbb{Q}[[x]].$$
Equality (\ref{2}) implies
\begin{equation}
\label{3}
{\frak c}'+(1\otimes {\frak g})(1\otimes x)=
(\Delta{\frak g})({\frak F}(x\otimes 1,1\otimes x)),
\end{equation}
where ${\frak c}'$ is independent of $1\otimes x.$
The application of $\id\otimes \widetilde{\varepsilon}$
to relation (\ref{3}) yields
$$
(\id \otimes \widetilde{\varepsilon}){\frak c}'=(((\id \otimes
\varepsilon)\circ \Delta){\frak g})(x\otimes 1)=({\frak g}\otimes
1)(x\otimes 1),
$$
and the application of $\; \widetilde{\varepsilon}\otimes \id \;$
to relation (\ref{3}) yields
$$
(\widetilde{\varepsilon}\otimes \id){\frak c}'
+(1\otimes{\frak g})(1\otimes x)=(1\otimes{\frak g})(1\otimes x).
$$
Hence
\begin{equation}
\label{4}
(\Delta{\frak g})({\frak F}(x\otimes 1,1\otimes x))={\frak c}+
({\frak g}\otimes 1)(x\otimes 1)+(1\otimes {\frak g})(1\otimes x),
\end{equation}
where ${\frak c}'=({\frak g}\otimes 1)(x\otimes 1)+{\frak c},\;
{\frak c}\in H_\mathbb{Q}{\mathop{\widehat{\otimes}}
\limits_{R_\mathbb{Q}}}H_\mathbb{Q}$ and $(\id \otimes \varepsilon)
{\frak c}=0=(\varepsilon \otimes \id){\frak c}.$

To complete the proof we must check the condition (i)
of the previous Lemma. For this purpose we apply
$\id \otimes \widetilde{\Delta}\; \hbox{and}\; \widetilde{\Delta}
\otimes \id$ to equation ({\ref{4}}).
We have
$$
(((\id \otimes \Delta)\circ \Delta){\frak g})
((\id \otimes \Delta){\frak F}(x
\otimes 1\otimes 1,1\otimes {\frak F}(x\otimes 1,1\otimes x)))=
$$
$$
(\id \otimes \Delta){\frak c}+{\frak g}(x)\otimes 1\otimes 1+
1\otimes (\Delta {\frak g})({\frak F}(x\otimes 1,1\otimes x))=
$$
$$
(((\Delta \otimes \id)\circ \Delta){\frak g})
((\Delta \otimes \id){\frak F}({\frak F}(x\otimes 1,1\otimes x)
\otimes 1,1\otimes 1 \otimes x))=
$$
$$
(\Delta \otimes \id){\frak c}+(\Delta {\frak g})
({\frak F}(x\otimes 1,1\otimes x))\otimes 1+1\otimes 1\otimes
{\frak g}(x),
$$
i. e.
$$(\id \otimes \Delta){\frak c}+{\frak g}(x)\otimes 1
\otimes 1+1\otimes {\frak c}+1\otimes {\frak g}(x)\otimes 1+
1\otimes 1\otimes {\frak g}(x)=
$$
$$(\Delta \otimes \id ){\frak c}+{\frak c}\otimes 1+
{\frak g}(x)\otimes 1\otimes 1+1\otimes {\frak g}(x)
\otimes 1+1\otimes 1\otimes {\frak g}(x).
$$
This completes the proof.\; $\square$
\begin{lemma}
``Linear`` formal groups ${\frak c}_i+x\otimes 1+1\otimes x,\; i=1,2$
over $H_\mathbb{Q}$ are isomorphic if and only if
cohomology classes $[{\frak c}_1]
\; \hbox{and}\; [{\frak c}_2]$ are equal.
\end{lemma}
{\raggedright {\it Proof}.}\quad
Suppose $[{\frak c}_1]=[{\frak c}_2];$ then there exists
$\lambda \in H_\mathbb{Q}$ such that $\varepsilon(\lambda)=0$
and ${\frak c}_2-{\frak c}_1=\Delta \lambda -\lambda \otimes 1-
1\otimes \lambda.$ Hence $\Delta \lambda +{\frak c}_1+x\otimes 1+1\otimes x=
{\frak c}_2+(\lambda+x)\otimes 1+1\otimes(\lambda +x).$
This shows that ${\frak g}(x)=\lambda +x$ is an isomorphism
from ${\frak c}_1+x\otimes 1+1\otimes x$ to
${\frak c}_2+x\otimes 1+1\otimes x.$
The proof of the converse statement
is clear.\quad $\square$\\

We may obtain more precise result than in the previous
Proposition.
\begin{corollary}
A formal group ${\frak F}(x\otimes 1,1\otimes x)$
over $H_\mathbb{Q}$ is isomorphic to the ``trivial`` group
$x\otimes 1+1\otimes x$ if and only if the 2-cocycle ${\frak c}$
is a coboundary.
\end{corollary}
{\raggedright {\it Proof}.}\quad
Let $$(\Delta{\frak g})({\frak F}(x\otimes 1,1\otimes x))={\frak c}+
{\frak g}(x)\otimes 1+1\otimes {\frak g}(x).$$
Let ${\frak c}=\lambda \otimes 1-\Delta \lambda +1\otimes \lambda.$
Let us consider the isomorphism ${\frak h}(x)=\lambda +{\frak g}(x).$
We have $$(\Delta{\frak h})({\frak F}(x\otimes 1,1\otimes x))=
\Delta \lambda +(\Delta{\frak g})({\frak F}(x\otimes 1,1\otimes x))=$$
$$\Delta \lambda+{\frak c}+{\frak g}(x)\otimes 1+1\otimes {\frak g}(x)
=\Delta \lambda -\lambda \otimes 1-1\otimes \lambda +{\frak c}+
{\frak h}(x)\otimes 1+1\otimes {\frak h}(x)=$$
$${\frak h}(x)\otimes 1+1\otimes {\frak h}(x).\quad \square$$

\begin{remark}
Note that this proof generalizes the standard proof of
the analogous result for formal groups over rings
(see \cite{Honda}).
\end{remark}
\begin{remark}
Note that in the proof we assign for any formal group
${\frak F}(x\otimes 1,1\otimes x)$ over $H$ some 2-cocycle ${\frak c}$
in the cobar complex of the coalgebra $H_\mathbb{Q}.$
\end{remark}
\begin{remark}
Note that $(\varepsilon{\frak g})(x)\in R_\mathbb{Q}[[x]]$
is the logarithm of the formal group $((\varepsilon \otimes
\varepsilon){\frak F})(x\otimes 1,1\otimes x)=F(x\otimes 1,
1\otimes x)\in R[[x\otimes 1,1\otimes x]]$
over ring $R.$
\end{remark}
\begin{remark}
Since ${\frak g}(x)=b_0x+b_1x^2+\ldots\quad \hbox{and}\; \varepsilon(b_0)=1,$
there exists the series
$(\Delta{\frak g})^{-1}(x)=(\Delta({\frak g}^{-1}))(x))\in
H_\mathbb{Q}{\mathop{\widehat{\otimes}}
\limits_{R_\mathbb{Q}}}H_\mathbb{Q}[[x]].$
Using (\ref{4}), we get
$$
{\frak F}(x\otimes 1,1\otimes x)=(\Delta{\frak g})^{-1}
({\frak c}+{\frak g}(x)\otimes 1+1\otimes {\frak g}(x)).
$$
\end{remark}

\section{The tensor category, connected with cobordisms rings of
FBSP}

Below
we study some properties of category, connected with
cobordism rings of FBSP. In particular, we shall show that
it is the tensor category.

Let us consider the series ${\frak G}(x,y)\in H[[x,y]]=\Omega^*_U
(\widetilde{Gr}),$ where $H=\Omega_U^*(Gr).$
Recall that it corresponds to the direct limit $\kappa$
of the maps $\kappa_{k,l}\colon \widetilde{Gr}_{k,kl}\rightarrow
\mathbb{C}P^{kl-1},$ where $\widetilde{Gr}_{k,kl}$
is the canonical FBSP over $Gr_{k,kl}$ ($(k,l)=1$).
Previously some properties of ${\frak G}(x,y)$ were
studied. In particular, it was shown that
$$
(\varepsilon{\frak G})(x,y)=F(x,y),
$$
where $\varepsilon \colon H\rightarrow R=\Omega^*_U(\pt)$
is the counit of the Hopf algebra $H$ and $F(x,y)\in R[[x,y]]$
is the formal group of geometric cobordisms.

Let $\varphi_{k,l}$ be the map
$$
\kappa_{k,l}\times \id_{\widetilde{Gr}_{k,kl}}\colon
\widetilde{Gr}_{k,kl}\rightarrow \mathbb{C}P^{kl-1}\times
\widetilde{Gr}_{k,kl}.
$$
The commutativity of the following diagram
\begin{equation}
\begin{array}{ccc}
\widetilde{Gr}_{k,kl} & \stackrel{\varphi_{k,l}}{\rightarrow} &
\mathbb{C}P^{kl-1}\times \widetilde{Gr}_{k,kl} \\
{\scriptstyle \varphi_{k,l}}\downarrow \quad \, && \qquad \quad \downarrow
{\scriptstyle \id_{\mathbb{C}P}\times \varphi_{k,l}} \\
\mathbb{C}P^{kl-1}\times \widetilde{Gr}_{k,kl} &
\stackrel{\diag_{\mathbb{C}P}\times \id_{\widetilde{Gr}}}{\rightarrow}
 & \mathbb{C}P^{kl-1}\times \mathbb{C}P^{kl-1}\times \widetilde{Gr}_{k,kl} \\
\end{array}
\end{equation}
allows us to define on the algebra $H[[x,y]]$ the structure of
$R[[z]]=\Omega^*_U(\mathbb{C}P^\infty)$-module
such that $z$ acts as the multiplication by ${\frak G}(x,y).$
Let us denote this $R[[z]]$-module by $(H[[x,y]];\: {\frak G}(x,y)).$

Let us consider $R[[z]]=\Omega^*_U(\mathbb{C}P^\infty)$ as a
Hopf algebra. Recall that $\Delta_{R[[z]]}(z)=F(z\otimes 1,1\otimes z).$

\begin{proposition}
$H[[x,y]]$ is the module coalgebra over $R[[z]],$
i. e. $R[[z]]{\mathop{\widehat{\otimes}}\limits_R}H[[x,y]]\rightarrow
H[[x,y]]$ is the homomorphism of coalgebras.
\end{proposition}
{\raggedright {\it Proof}.}\quad The proof follows from the
following commutative diagram
($(km,ln)=~1$):
\begin{equation}
\nonumber
\begin{array}{ccc}
\diagram
& {\scriptstyle \mathbb{C}P^{kl-1}\times \widetilde{Gr}_{k,kl}\times
\mathbb{C}P^{mn-1}\times \widetilde{Gr}_{m,mn}} \dlto & \\
{\scriptstyle\mathbb{C}P^{kl-1}\times \mathbb{C}P^{mn-1}\times
\widetilde{Gr}_{k,kl}\times \widetilde{Gr}_{m,mn}} \dto &&
{\scriptstyle\widetilde{Gr}_{k,kl}\times
\widetilde{Gr}_{m,mn}} \ulto_{\; \varphi_{k,l}\times \varphi_{m,n}}\dto & \\
{\scriptstyle\mathbb{C}P^{klmn-1}\times \widetilde{Gr}_{km,klmn}} &&
\quad{\scriptstyle\widetilde{Gr}_{km,klmn}}.\quad \square \llto_{\varphi_{km,ln}}
\enddiagram
\end{array}
\end{equation}

Let us consider the next commutative diagram ($(km,ln)=1$):
\begin{equation}
\begin{array}{ccc}
G\widetilde{r_{k,kl}\times Gr}_{m,mn} & \stackrel{\psi_{kl,mn}}
{\rightarrow} & \widetilde{Gr}_{km,klmn} \\
\downarrow && \downarrow \\
Gr_{k,kl}\times Gr_{m,mn} & \stackrel{\phi_{kl,mn}}{\rightarrow} &
\; Gr_{km,klmn}\; ,
\end{array}
\end{equation}
where $G\widetilde{r_{k,kl}\times Gr}_{m,mn}$ is the FBSP
over $Gr_{k,kl}\times Gr_{m,mn},$ induced by the map $\phi_{kl,mn}.$
Clearly that the bundle $G\widetilde{r_{k,kl}\times Gr}_{m,mn}$
(with fiber $\mathbb{C}P^{km-1}\times \mathbb{C}P^{ln-1}$)
is (``external``) Segre's product of the canonical FBSP over $Gr_{k,kl}$
and $Gr_{m,mn}.$ By definition, put
$$
\widetilde{\; Gr\times
Gr}=\varinjlim_{(km,ln)=1}G\widetilde{r_{k,kl}\times Gr}_{m,mn}\; ,
$$
$$
\psi=\varinjlim_{(km,ln)=1}\psi_{km,ln}\colon \quad \widetilde{Gr\times Gr}
\rightarrow \widetilde{Gr}.
$$

We have the homomorphism of $R[[z]]$-modules
$$
\Psi \colon (H[[x,y]];\: {\frak G}(x,y))\rightarrow (H
{\mathop{\widehat{\otimes}}\limits_R}H[[x,y]];\: (\Delta{\frak
G})(x,y))\: ,
$$
defined by the fiber map $\psi$
(recall that $\Delta$ is the comultiplication in the Hopf algebra
$H=\Omega^*_U(Gr)$).
Clearly that the restriction $\Psi \mid_H$ coincides with $\Delta.$

Let ${\cal P}^{k-1}{\mathop{\times}\limits_X}{\cal Q}^{l-1}$
be a FBSP over a finite $CW$-complex $X$
with fiber $\mathbb{C}P^{k-1}\times \mathbb{C}P^{l-1}.$
Recall that if $k$ and $l$ are sufficiently large then
there exist a classifying map $f_{k,l}$
and the corresponding fiber map
\begin{equation}
\begin{array}{ccc}
{\cal P}^{k-1}{\mathop{\times}\limits_X}{\cal Q}^{l-1} &
\rightarrow & \widetilde{Gr}_{k,kl} \\
\downarrow && \downarrow \\
X & \stackrel{f_{k,l}}{\rightarrow} & Gr_{k,kl}
\end{array}
\end{equation}
which are unique up to homotopy and up to fiber homotopy respectively.
Let ${\cal P}^{km-1}{\mathop{\times}\limits_X}{\cal Q}^{ln-1},\;
(km,ln)=1$
be Segre's product of
${\cal P}^{k-1}{\mathop{\times}\limits_X}{\cal Q}^{l-1}$
with the trivial FBSP $X\times \mathbb{C}P^{m-1}\times \mathbb{C}P^{n-1}.$
Let us pass to the direct limit
$${\cal P}{\mathop{\times}\limits_X}{\cal Q}=
\varinjlim_i({\cal P}^{km_i-1}{\mathop{\times}\limits_X}
{\cal Q}^{ln_i-1}),$$ where $(km_i,ln_i)=1,$
$m_i\mid m_{i+1},\: n_i\mid n_{i+1},\; m_i\, ,n_i\rightarrow
\infty,$ as $i\rightarrow \infty.$
The stable equivalence class of FBSP
${\cal P}^{k-1}{\mathop{\times}\limits_X}{\cal Q}^{l-1}$
may be unique restored by the direct limit
${\cal P}{\mathop{\times}\limits_X}{\cal Q}.$
We have also a classifying map $f=\varinjlim_{(k,l)=1}f_{k,l}$
and the corresponding fiber map
\begin{equation}
\begin{array}{ccc}
{\cal P}{\mathop{\times}\limits_X}{\cal Q} &
\rightarrow & \widetilde{Gr} \\
\downarrow && \downarrow \\
X & \stackrel{f}{\rightarrow} & Gr\; .
\end{array}
\end{equation}

Let us define the category $\frak{FBSP}_{finite}$
by the following way.
\begin{itemize}
\item[(i)] $\Ob (\frak{FBSP}_{finite})$ is the class of direct limits
${\cal P}{\mathop{\times}\limits_X}{\cal Q}$ of FBSP over finite
$CW$-complexes $X$ (in other words, the class of stable equivalence
classes of FBSP);
\item[(ii)] $\Mor_{\frak{FBSP}_{finite}}({\cal P}{\mathop{\times}\limits_X}
{\cal Q},\: {\cal P}'{\mathop{\times}\limits_Y}{\cal Q}')$
is the set of fiber maps
\begin{equation}
\begin{array}{ccc}
{\cal P}{\mathop{\times}\limits_X}{\cal Q} &
\rightarrow & {\cal P}'{\mathop{\times}\limits_Y}{\cal Q}' \\
\downarrow && \downarrow \\
X & \rightarrow & Y
\end{array}
\end{equation}
such that its restrictions to any fiber
$(\cong \mathbb{C}P^\infty\times \mathbb{C}P^\infty)$ are isomorphisms.
\end{itemize}

Applying the functor of unitary cobordisms $\Omega^*_U$
to an object ${\cal P}{\mathop{\times}\limits_X}{\cal Q}
\in \Ob(\frak{FBSP}_{finite}),$
we get the $R[[z]]$-module $(A[[x,y]];\:(f^*{\frak G})(x,y))\in
\Ob(\Omega^*_U(\frak{FBSP}_{finite})),$ where $A=\Omega^*_U(X)$ and $f
\colon X\rightarrow Gr$ is a classifying map for
${\cal P}{\mathop{\times}\limits_X}{\cal Q}.$
It is clear that $((\varepsilon_A\circ f^*){\frak G})(x,y)=F(x,y),$
where $\varepsilon_A \colon A\rightarrow R$ is the homomorphism,
induced by an embedding of a point $\pt \hookrightarrow X.$
In other words, for any object in the category
$\Omega^*_U(\frak{FBSP}_{finite})$
there exists the canonical morphism $(A[[x,y]];\:(f^*{\frak G})(x,y))
\rightarrow (R[[x,y]];\: F(x,y)).$

Hence there exist the initial object $(H[[x,y]];\: {\frak G}(x,y))$
and the final object $(R[[x,y]];\: F(x,y))$ in the category
$\Omega^*_U(\frak{FBSP}).$

Let's consider a pair $(A[[x,y]];\:(f^*{\frak G})(x,y)),\;
(B[[x,y]];\:(g^*{\frak G})(x,y))\in \Ob(\Omega^*_U(\frak{FBSP}_{finite})),$
where $(B[[x,y]];\:(g^*{\frak G})(x,y))=
\Omega^*_U({\cal P}'{\mathop{\times}\limits_Y}{\cal Q}').$
Let's define their ``tensor product`` as the object
$((A{\mathop{\otimes}\limits_R}B)[[x,y]];\:(((f^*\otimes g^*)\circ
\Delta){\frak G})(x,y))\in \Ob(\Omega^*_U(\frak{FBSP}_{finite}))$
(recall that $\Delta \colon H\rightarrow
H{\mathop{\widehat{\otimes}}\limits_R}H$ is the comultiplication
in the Hopf algebra $H$).

\begin{proposition}
The category $\Omega^*_U(\frak{FBSP}_{finite})$ is the tensor category
with the just defined tensor product and the unit
$\; \I =(R[[x,y]];\: F(x,y)).$
\end{proposition}
{\raggedright {\it Proof}.}\quad The proof is trivial.
For example, the associativity axiom follows from the
identity $(((\Delta \otimes \id_H)\circ \Delta){\frak G})(x,y)=
(((\id_H\otimes \Delta)\circ \Delta){\frak G})(x,y)$
which follows from the next commutative diagram ($(kmt,lnu)=1$):
\begin{equation}
\begin{array}{ccc}
G\widetilde{r_{km,klmn}\times Gr}_{t,tu} &
\rightarrow & \widetilde{Gr}_{kmt,klmntu} \\
\uparrow && \uparrow \\
Gr_{k,kl}\widetilde{\times Gr_{m,mn}}
\times Gr_{t,tu} &
\rightarrow & G\widetilde{r_{k,kl}\times G}r_{mt,mntu},
\end{array}
\end{equation}
where $Gr_{k,kl}\widetilde{\times Gr_{m,mn}}
\times Gr_{t,tu}$
is external Segre's product of the canonical FBSP over
$Gr_{k,kl},\; Gr_{m,mn}$ and $Gr_{t,tu}$ (it
is the bundle over $Gr_{k,kl}\times Gr_{m,mn}\times Gr_{t,tu}$ with
fiber $\mathbb{C}P^{kmt-1}\times \mathbb{C}P^{lnu-1}$).\:$\square$

\smallskip

Note that there exist the canonical homomorphisms $p_1,\; p_2$:
$$
\diagram
& \scriptstyle{((A{\mathop{\otimes}\limits_R}B)[[x,y]];\:
(((f^*\otimes g^*)\circ
\Delta){\frak G})(x,y))} \dlto_{\scriptscriptstyle{p_1}}
\drto^{\scriptscriptstyle{p_2}} &\\
\scriptstyle{(A[[x,y]];\:(f^*{\frak G})(x,y))} &&
\scriptstyle{(B[[x,y]];\:(g^*{\frak
G})(x,y))}\; ,\\
\enddiagram
$$
such that $p_1\mid_{A{\mathop{\otimes}\limits_R}B}=\id_A\otimes
\varepsilon_B,\quad
p_2\mid_{A{\mathop{\otimes}\limits_R}B}=\varepsilon_A
\otimes \id_B.$

\smallskip

The author is grateful to professor E. V. Troitsky
for constant attention to this work, and to professors
V. M. Manuilov and A. S. Mishchenko for useful discussions.

\end{document}